\documentclass[11pt]{amsart}

\usepackage[colorlinks=true,citecolor=cyan,backref=page]{hyperref}

\usepackage[square,sort&compress,comma,numbers]{natbib}
\usepackage{nicefrac,xcolor,upref}

\usepackage{amscd,amsthm,amsmath,amssymb,amsfonts}

\usepackage{mathrsfs}

\newtheorem{theorem}{Theorem}[section]
\newtheorem{lemma}[theorem]{Lemma} 
 
\newtheorem{definition}[theorem]{Definition}

\numberwithin{equation}{section}

\def\hw{{\widehat w}} 
\def\hv{{\widehat v}}

\def\Z{{\mathbb Z}}  
\def\R{{\mathbb R}}

\def\eps{{\varepsilon}}

\def\cR{{\mathcal R}}

\def\tx{{\widetilde x}}

\def\bfx{{\bf x}}

\def\beq{\begin{equation}}
\def\eeq{\end{equation}}

\def\ice{{\rm ice}}
\def\hice{{\widehat \ice}}

\def\beq{\begin{equation}}
\def\eeq{\end{equation}}

\def\hr{{\widehat r}}

\def\Per{{\rm Per}}

\begin{document}

\title{Remarks on uniform recurrence properties for beta-transformation}

\author{Yann Bugeaud}
\address{I.R.M.A., UMR 7501, Universit\'e de Strasbourg
et CNRS, 7 rue Ren\'e Descartes, 67084 Strasbourg Cedex, France}
\address{Institut universitaire de France}
\email{bugeaud@math.unistra.fr}

\begin{abstract}
We complement the recent paper of Zheng and Wu 
[{\it Uniform recurrence properties for beta-transformation}, 
Nonlinearity 33 (2020), 4590--4612], where the authors study, from the 
metrical point of view, the uniform recurrence properties of the orbit of a point 
under the $\beta$-transformation to the point itself. 
\end{abstract}


\maketitle

\section{Introduction}\label{sec:1}

Let $\beta > 1$ be a real number. The \emph{$\beta$-transformation} $T_\beta$ on $[0,1)$ is defined by
\begin{align*}
T_{\beta}(x) = \beta x\mod 1.
\end{align*}
In a recent paper, Zheng and Wu \cite{ZhWu20} have introduced the following exponents of approximation $r_\beta$ 
and $\hr_\beta$. 

\begin{definition}
Let $\beta > 1$ be a real number. 
For $x$ in $[0,1)$, set
$$
r_{\beta}(x):=\sup \{r \ge 0 : |T_{\beta}^nx-x|<(\beta^n)^{-r}\ {\rm for \ infinitely\ many \ integers \ }  n \ge 1\}
$$ 
and 
\begin{align*}
\widehat{r}_\beta(x) & :=\sup\{ \hr \ge 0 :\ {\rm for\ every\ sufficiently\ large\ integer}\ N, \\
& {\rm there\ exists\ }n \ {\rm such\ that\ } 1 \le n \le N \ {\rm and} \ |T_{\beta}^nx-x|<(\beta^N)^{-\widehat{r}}\}.
\end{align*}
\end{definition}

Said differently, $r_{\beta}(x)$ is the supremum of the real numbers $r$ for which the inequality
$$
|T_{\beta}^nx-x|<(\beta^n)^{-r}
$$ 
has infinitely many solutions in positive integers $n$, while $\widehat{r}_{\beta}(x)$ is 
the supremum of the real numbers $\widehat{r}$ such that, for every large enough integer $N$, 
the inequality
$$
|T_{\beta}^nx-x|<(\beta^N)^{-\widehat{r}}
$$ 
has a solution $n$ with $1\leq n\leq N$. 
The exponents $r_\beta$ and $\hr_\beta$ should be compared with the exponents $v_\beta$ and 
$\hv_\beta$ introduced in \cite{BuLi16}, which are defined analogously, but with $|T_{\beta}^nx-x|$ replaced 
by $|T_{\beta}^n x|$. Roughly speaking, $v_\beta (x)$ and $\hv_\beta (x)$ 
measure the sizes of large blocks of $0$ (or of the 
digit $\lceil \beta \rceil - 1$) in the $\beta$-expansion of $x$,
while $r_\beta (x)$ 
and $\hr_\beta (x)$ measure repetitions of prefixes in the $\beta$-expansion of $x$. 
Throughout, we let $\lceil \cdot \rceil$ and $\lfloor \cdot \rfloor$ denote the upper integral part and the 
lower integral part, respectively. 

For $r$ and $\hr$ with $0\leq r, \hr \le +\infty$, set
$$
R_\beta(\widehat{r},r):=\left\{x \in [0,1):\widehat{r}_\beta(x)=\widehat{r},\ r_\beta(x)=r\right\}.
$$
Observe that $R_\beta (+ \infty, + \infty)$ is equal to the (countable) set $\Per (\beta)$ of real numbers $x$ in $[0, 1)$ whose 
$\beta$-expansion is purely periodic.  
Zheng and Wu \cite[Theorem 1.2]{ZhWu20} 
have determined the Hausdorff dimension of the sets $R_\beta(\widehat{r},r)$ when $\hr \le r / (1 +r)$. They also claimed that 
the sets $R_\beta(\widehat{r},r)$ are countable when $\hr > r / (1 +r)$.
In the present note, we explain why this is not the case and correct their result. 
We establish the following theorems.

\begin{theorem}   \label{sturm}
Let $\beta > 1$ be a real number. 
There exists an uncountable set of real numbers $\cR$
such that, for every $r$ in $\cR$, there are uncountably many real numbers $x$ with 
$\hr_\beta (x) = 1$ and $r_\beta (x) = r$. 
\end{theorem}

Theorem \ref{sturm} is a consequence of classical results on Sturmian sequences recalled in Section \ref{secSturm}. 
It shows that the set $R_\beta (1, r)$ is uncountable for uncountably many values of $r$ and that, consequently, 
the second assertion of \cite[Theorem 1.2]{ZhWu20} is incorrect. 
To fully characterize the triples $(\beta, \hr, r)$ in $\R_{> 1} \times \R_{>0} \times \R_{>0}$ with 
$\hr > r / (1 +r)$  such that  the set $R_\beta(\widehat{r},r)$ is uncountable is a very difficult question, even when $\beta$ is an integer.

\begin{theorem}   \label{bound}
For every real number $\beta > 1$ and every real number $x$ in $[0, 1)$
whose $\beta$-expansion is not purely periodic, we have $\hr_\beta (x) \le 1$. 
\end{theorem}

Recalling that $\hr_\beta (x)$ is infinite for every real number $x$ in $[0, 1)$ whose 
$\beta$-expansion is purely periodic, Theorem \ref{bound} shows that, for every real number 
$\hr > 1$, the set $\{x \in [0, 1) : \hr_\beta (x) = \hr \}$ is empty. 
Theorem \ref{bound} is well-known when $\beta$ is an integer; see Section \ref{secSturm}. 
There is, however, a small difficulty to overcome 
to treat the general case.

\begin{theorem}   \label{thdim}
For every real number $\beta > 1$, the set of real numbers $x$ in $[0, 1)$ 
such that $\hr_\beta (x) > r_\beta (x) / (1 + r_\beta (x))$ has 
Hausdorff dimension $0$. 
\end{theorem}

It follows from Theorem \ref{thdim} that all the metrical results of \cite{ZhWu20} are correct. 
In particular,
for $\beta > 1$ and $\hr$ in $[0, 1]$, we have
$$
\dim \{x\in[0,1):\widehat{r}_\beta(x)\geq\widehat{r}\}
=\dim \{x\in[0,1):\widehat{r}_\beta(x)=\widehat{r}\}=\left(\frac{1-\widehat{r}}{1+\widehat{r}}\right)^2, 
$$
where $\dim$ stands for the Hausdorff dimension.

The present paper is organized as follows. We discuss $\beta$-expansions from a combinatorial point of view in 
Section \ref{sec:3}, mostly following \cite{ZhWu20}. 
We recall classical results on Sturmian sequences in Section \ref{secSturm} and use them to establish 
Theorem \ref{sturm}. 
The last two sections are devoted to the proofs of Theorems \ref{bound} and \ref{thdim}, respectively.

\section{Diophantine approximation and combinatorics on words}  \label{sec:3}

Every real number $x$ in $[0, 1)$ can be uniquely expanded as
$$
x = \frac{\varepsilon_1(x,\beta)}{\beta} +\cdots+\frac{\varepsilon_n(x,\beta)}{\beta^n}+\cdots, 
$$
where $\varepsilon_n(x,\beta)=\lfloor\beta T_{\beta}^{n-1}(x)\rfloor$ 
for $n \ge 1$.
The integer $\varepsilon_n(x,\beta)$ is called the \emph{$n$-th digit of $x$}. We call the sequence 
$\varepsilon(x,\beta):= ((\varepsilon_ n(x,\beta))_{n \ge 1}$ the \emph{$\beta$-expansion of $x$}. 
When there is no ambiguity, we simply write $\eps_n$ for $\eps_n (x, \beta)$. 

We now extend the definition of the $\beta$-transformation to $x=1$. 
Let $T_\beta(1)=\beta-\lfloor \beta\rfloor$. We have 
$$
1=\frac{\varepsilon_1(1,\beta)}{\beta}+\cdots+\frac{\varepsilon_n(1,\beta)}{\beta^n}+\cdots,
$$ 
where $\varepsilon_n(1,\beta)=\lfloor\beta T_\beta^n(1)\rfloor.$ 
If the $\beta$-expansion of $1$ is finite, that is, if there is an integer $m\geq 1$ 
such that $\varepsilon_m(1,\beta)> 0$ and  
$\varepsilon_k(1,\beta)=0$ for all $k> m$, 
then $\beta$ is called a \emph{simple Parry number}. In this case, set 
$$
\varepsilon^\ast(\beta):=(\varepsilon_1^\ast(\beta),\varepsilon_2^\ast(\beta),\ldots)
=(\varepsilon_1(1,\beta),\varepsilon_2(1,\beta),\ldots, \varepsilon_m(1,\beta)-1)^\infty , 
$$ 
where $\omega^\infty=(\omega,\omega,\ldots)$. If the $\beta$-expansion of $1$ is not finite, 
set $\varepsilon^\ast(\beta)=\varepsilon(1,\beta)$. In both cases, we have 
\beq  \label{infbeta}
1=\frac{\varepsilon_1^\ast(\beta)}{\beta}+\cdots+\frac{\varepsilon_n^\ast(\beta)}{\beta^n}+\cdots.
\eeq
The sequence $\varepsilon^\ast(\beta)$ is called \emph{the infinite $\beta$-expansion of $1$}.

As in \cite{ZhWu20}, for any $x$ in $[0, 1)$ whose $\beta$-expansion is not periodic and such that $r_\beta (x) > 0$, 
there exist two integer sequences   $(n_k)_{k\ge 1}$ and $(m_{k})_{k\ge 1}$ such that
$(m_k - n_k)_{k \ge 1}$ is increasing, 
\begin{align*}
r_\beta(x)=\limsup_{k\rightarrow + \infty}\frac{m_{k}-n_{k}}{n_{k}}
\end{align*}
and
\begin{equation}\label{inf}
\widehat{r}_\beta(x)=\liminf_{k\rightarrow + \infty}\frac{m_{k}-n_{k}}{n_{{k}+1}}.
\end{equation}
To see this, let us set
$$
n'_1=\min\{n\geq1:\varepsilon_{n+1}=\varepsilon_1\},\ \ m'_1=\max\{n\geq n_1: |T_\beta^{n'_1} x-x|< \beta^{-(n-n'_1)}\}.
$$ 
Suppose that, for some $k\geq 1$, the integers $n'_k$ and $m'_k$ have been defined. Then, set 
\begin{align*}
n'_{k+1} & = \min\{n\geq n'_k:\varepsilon_{n+1}=\varepsilon_1\}, \\
m'_{k+1} & =\max\{n\geq n'_{k+1}: |T_\beta^{n'_k} x-x|< \beta^{-(n-n'_k)}\}.
\end{align*}
Since $r_\beta(x)>0$, the digit $\varepsilon_1$ occurs infinitely often in the word $\varepsilon(x,\beta)$, thus
$n'_k$ is well 
defined. Since $\varepsilon(x,\beta)$ is not periodic, $m'_k$ is well defined. By the definitions of $n'_k$ and $m'_k$,  
for all $k\geq1$, we have 
$$
\beta^{-(m'_k-n'_k)-1}\leq |T_\beta^{n'_k} x-x|< \beta^{-(m'_k-n'_k)}.
$$
Now, put $i_1=1$. Assume that $i_k$ has been defined. Let 
$$
i_{k+1}=\min\{i>i_k: m'_i -n'_i>m'_{i_k}-n'_{i_k}\}.
$$
For $k \ge 1$, we put $m_k = m'_{i_k}$ and $n_k = n'_{i_k}$. In particular, the sequence $(m_k - n_k)_{k \ge 1}$ 
is increasing.

The exponent $r_\beta$ is close to the initial critical exponent (ice), which has been introduced 
by Berth\'e, Holton, and Zamboni \cite{BeHoZa06}. 
For an infinite word $\bfx$ over a finite alphabet, 
we let $\ice(\bfx)$ (resp., $\hice(\bfx)$) denote the supremum of the real numbers $\rho$ 
such that there are arbitrarily large integers $N$ (resp., for all sufficiently large $N$), 
there exists $n$ such that $1 \le n \le N$ and the prefix of $\bfx$ 
of length $\lfloor \rho n \rfloor$ is a power of the prefix of $\bfx$ of length $n$. 
Observe that $\ice(\bfx) = \hice(\bfx) = + \infty$ if and only if $\bfx$ is purely periodic.

It is easy to see that, 
for any $\beta > 1$ and any $x$ in $[0, 1)$, we have
$$
r_\beta (x) \ge \ice(\eps(x, \beta)) - 1, \quad \hr_\beta (x) \ge \hice(\eps(x, \beta)) - 1.
$$
At first sight, it may seem that both inequalities are equalities. However, this is not the case (at least for the first one), 
because there are two admissible ways to write $1$ as a sum of powers of $1/\beta$, namely 
$1 + \sum_{n \ge 1} 0 / \beta^n$ and the infinite $\beta$-expansion of $1$ defined by \eqref{infbeta}. 
The situation is described in the following lemma, extracted from \cite{ZhWu20}.

\begin{lemma}  \label{tkmk}
Let $\beta > 1$ be a real number. 
Let $x$ be in $[0,1)$ with 
$r_\beta(x)>0$ and whose $\beta$-expansion is not periodic. 
Let $(m_k)_{k \ge 1}$ and $(n_k)_{k \ge 1}$ 
be as above. Set
\begin{align*}
t_k=\max\{n>n_k:(\varepsilon_{n_k+1}\ldots,\varepsilon_n) =(\varepsilon_1,\ldots,\varepsilon_{n-n_k})\}.
\end{align*}
Then $t_k \leq m_k + 1$. If $t_k \ge m_k$, then we have
\begin{align*}
\left(\varepsilon_1,\ldots,\varepsilon_{m_{k}}\right)=\left(\varepsilon_1,\ldots,\varepsilon_{n_{k}},\varepsilon_1,\ldots, \varepsilon_{m_{k}-n_{k}}\right).
\end{align*} 
If $t_k<m_k$, then we have 
\begin{align*}
\begin{split} \label{1}
(\varepsilon_1,\ldots,\varepsilon_{m_k})
=\bigl(\varepsilon_1, & \ldots,\varepsilon_{n_k},\varepsilon_1,\ldots,  \varepsilon_{t_k-n_k}, \\ 
& \varepsilon_{t_k-n_k+1}-1, 
\varepsilon_1^\ast(\beta), \ldots, \varepsilon_{m_k-t_k-1}^\ast(\beta)\bigr) 
\end{split}
\end{align*} 
or 
\begin{align*}
(\varepsilon_1,\ldots,\varepsilon_{m_k})
=(\varepsilon_1,\ldots,\varepsilon_{n_k},\varepsilon_1,\ldots,\varepsilon_{t_k - n_k},
\varepsilon_{t_k-n_k+1}+1, 0^{m_k-t_k-1}).
\end{align*}
\end{lemma}

Let us observe that, if $t_k<m_k$ in Lemma \ref{tkmk}, then the first case occurs when $\eps_{t_k + 1} > 0$ and 
$\eps_{t_k + 2} = \ldots = \eps_{m_k} = 0$, while the second case occurs when $\eps_{t_k + 1} < \lceil \beta \rceil - 1$ and 
$(\eps_{t_k + 2},  \ldots ,  \eps_{m_k} ) = ( \eps_1^* (\beta), \ldots , \eps^*_{m_k - t_k -1} (\beta) )$. 

Let us note that, when the sequence $(m_k - t_k)_{k \ge 1}$ is bounded, 
then we have $r_\beta (x) = \ice ( \eps(x, \beta))$ and  $\hr_\beta (x) = \hice ( \eps(x, \beta))$.

\section{Proof of Theorem \ref{sturm}}   \label{secSturm}

\subsection{Connection between the exponent $r_\beta$ and rational approximation, when 
$\beta$ is an integer}

Let $x$ be an irrational real number in $[0, 1)$. 
We denote by $w_1 (x)$ (resp., $\hw_1 (x)$) the supremum of the 
real numbers $w_1$ such that, for arbitrarily large integers $Q$ 
(resp., for every sufficiently large integer $Q$), there exists $q$
with $1 \le q \le Q$ and $\| q x \| \le Q^{-w_1}$, where $\| \cdot \|$ is the distance to the nearest integer. 
While $w_1$ can take every value $\ge 2$, the function $\hw_1$ is constant equal to $1$ 
on the set of irrational numbers. This was proved by Khintchine \cite{Kh26}; see e.g. \cite{Bu16}.

Let $g \ge 2$ be an integer. Then, the transformation $T_g$ satisfies
$$
|T_g^n x - x | = \| g^n x - x \| = \| (g^n - 1) x \|, \quad n \ge 1.
$$
This shows that $r_g (x)$ measures the quality of approximation to $x$ by 
rational numbers whose denominator is of the form $g^n - 1$ for some positive integer $n$, that is, 
by rational numbers whose $g$-ary expansion is purely periodic. Consequently, we have 
\begin{equation}  \label{rgw1}
r_g (x) \le w_1 (x), \quad \hr_g (x) \le \hw_1 (x) = 1, \quad x \in [0, 1). 
\end{equation}
This establishes Theorem \ref{bound} when $\beta$ is an integer.

\subsection{Characteristic Sturmian sequences}   \label{subsecSturm}

We introduce a large class of real numbers for which we have equalities in \eqref{rgw1}. 
Let $(s_k)_{k\ge 1}$ be a sequence of positive integers and  
$\{a, b\}$ a two-letters alphabet. We define
inductively a sequence of finite words $(W_k)_{k \ge 0}$ over $\{a,b\}$ by
the formulas
$$
W_0 = b, \quad W_1 = b^{s_1 - 1} a,  \quad \hbox{and} \quad
W_{k+1} = W_k^{s_{k+1}} \, W_{k-1}, \quad k \ge 1.
$$
This sequence converges 
to the infinite word
$$
W_{\varphi} := \lim_{k \rightarrow \infty} \,
W_k = b^{s_1 - 1} a \ldots,
$$
which is usually called the {\it Sturmian characteristic word
of slope} 
$$
\varphi := [0; s_1, s_2, s_3, \ldots]
$$
constructed over the  alphabet $\{a, b\}$.
The characters $a$ and
$b$ will denote either the letters of the alphabet, or
distinct positive integers, according to the context. 
Let $g$ be an integer with $g \ge 2$. 
Let $\xi_{\varphi}$
denote the real number whose 
$g$-ary expansion is given by the letters of the infinite word $W_{\varphi}$ written over the two-letters alphabet 
$\{0, 1\}$, and set
$$
\sigma_{\varphi} =  \limsup_{k \to + \infty} \, [s_k; s_{k-1}, \dots, s_1].
$$
Then, we have
\beq     \label{rghat}
w_1 (\xi_\varphi) = r_g (\xi_\varphi) = \sigma_{\varphi} \quad \hbox{and} \quad  
\hw_1 (\xi_\varphi)  = \hr_g (\xi_\varphi) = 1.
\eeq
These results have been first established by B\"ohmer \cite{Bohm27}; see also \cite{AdDa77,BuLau23}. 
They imply that there exists an uncountable set of real numbers $\cR$
such that, for every $r$ in $\cR$, there are uncountably many real numbers $x$ with 
$\hr_g (x) = 1$ and $r_g (x) = r$. 
To see this, let $\psi : \Z_{\ge 1} \to \{3, 4\}$ and 
$$
\omega : \Z_{\ge 1} \setminus \{2^k - h : k \ge 2, 0 \le h < k\} \to \{1, 2\}
$$
be two arbitrary functions. Define the sequence $(s_k^{\psi, \omega})_{k \ge 1}$ by 
$$
s_{2^k}^{\psi, \omega} = 5, \quad s_{2^k - h}^{\psi, \omega} = \psi (h), \quad k \ge 2, \quad 1 \le h < k, 
$$
and
$$
s_\ell^{\psi, \omega} = \omega(\ell), \quad \ell \in \Z_{\ge 1} \setminus \{2^k - h : k \ge 2, 0 \le h < k\}.
$$
Then, $w_1 (\xi_\varphi) = [5; \psi(1), \psi(2), \ldots]$ does not depend on $\omega$. 
This completes the proof of Theorem \ref{sturm} when $\beta$ is an integer.  

\subsection{The Fibonacci sequence} 
Let us discuss a bit more the particular case of the Fibonacci sequence, obtained by taking 
for $(s_k)_{k \ge 1}$ the constant sequence equal to $1$. Then, the sequence $(W_k)_{k \ge 0}$ starts with
$$
b, \quad a, \quad ab, \quad aba, \quad abaab, \ldots 
$$
and
$$
W_\phi = abaababaabaababaababa \ldots ,
$$
where $\phi = (\sqrt{5} - 1) / 2$. 
Let $(F_k)_{k \ge 0}$ denote the Fibonacci sequence defined by $F_0 = 0$, $F_1 = 1$, and 
$F_{k+2} = F_{k+1} + F_k$ for $k \ge 0$. Then, we check that, for $k \ge 1$, the word $W_k$ is the 
prefix of $W_\phi$ of length $F_{k+1}$. Furthermore, the word $W_\phi$ starts with $W_k W_k W_{k-1}'$, where 
$W_{k-1}'$ is obtained from $W_{k-1}$ by permuting its last two letters. 
(Note that, for an arbitrary sequence $(s_k)_{k \ge 1}$ as above, we also have $W_k W_{k-1}' = W_{k-1} W_k$, and this is 
the key for the proof that $\hr_g (\xi_\varphi) = 1$.) 
The point here is that the length of $W_{k+1} = W_k W_{k-1}$ is (significantly) smaller than that of $W_k W_k W_{k-1}'$. 

This contrasts with the problem studied in \cite{BuLi16}, where the authors introduced the exponents 
$v_\beta$ and $\hv_\beta$ to measure the size of the blocks of the digit $0$ (or of the digit $\lceil \beta \rceil -1$) 
occurring near the beginning of the $\beta$-expansion of a given real number; see Section \ref{sec:1} above for their
precise definition. 
Obviously, two different such blocks do not overlap, giving after some 
computation that $\hv_\beta (x)$ is at most equal to $v_\beta (x)  / (1 + v_\beta (x))$, for every $x$ in $[0, 1)$.   
Here, the situation is different: we have just seen that $\hr_\beta (x)$ may exceed $r_\beta (x) / (1 + r_\beta (x))$. 
However, this happens only for very few $x$ in $[0, 1)$, by Theorem \ref{thdim}. 

\subsection{Completion of the proof of Theorem \ref{sturm}}

We keep the notation from Subsection \ref{subsecSturm}. 
Let $k$ be a positive integer. We work in the base $g = 2^k$. 
Consider the morphism $f_k$ from the set of words over $\{a, b\}$ to the set 
of words over $\{0, 1\}$ defined by $f_k(a) = 0^k$ and $f_k (b) = 0^{k-1} 1$. 
Let $\xi_{\varphi, k}$ be the real number whose $2^k$-ary expansion 
is given by the word $f_1 (W_\varphi)$.  Observe that its binary expansion 
is given by the word $V_{\varphi, k} := f_k (W_\varphi)$. 
Furthermore, the definition of $f_k$ shows that $\ice (W_\varphi) = \ice (V_{\varphi, k})$.

Let $\beta > 1$ be a real number which is not an integer. 
If $k$ is sufficiently large, then the word $V_{\varphi, k}$ and all its shifts are 
lexicographically smaller than $\eps^* (\beta)$. This implies, by \cite[Theorem 3]{Par60}, 
that the word $V_{\varphi, k}$ is admissible, that is, there exists a real number $x$ whose 
$\beta$-expansion is given by $V_{\varphi, k}$. Then, by \eqref{rghat}, we have
$$
r_{2^k} (\xi_{\varphi, k})  = r_{2} (\xi_{\varphi, k})  = r_\beta (x) = \sigma_\varphi \quad
{\rm and} \quad 
\hr_{2^k} (\xi_{\varphi, k}) = \hr_{2} (\xi_{\varphi, k}) = \hr_\beta (x) = 1. 
$$
This completes the proof of Theorem \ref{sturm} when $\beta$ is not an integer.

\section{Proof of Theorem \ref{bound}}  \label{proofTh2}

By \eqref{rghat}, we have $\hr_g (x) \le 1$, for every integer $g \ge 2$ and every irrational number in $[0, 1)$. 
This proves Theorem \ref{bound} when $\beta$ is an integer. 

Let $\beta > 1$ be a real number which is not an integer.  
Let $x$ be in $[0, 1)$  
and whose $\beta$-expansion is not 
ultimately periodic. 
Assume that $\hr_\beta (x) > 1$. 
By \eqref{inf}, there exists $\delta > 0$ such that, for $k$ large enough, say, for $k \ge k_0$, 
we have $n_{k+1} (1 +  \delta) \le  m_k - n_k$ and
$m_k > (2 + \delta) n_k$.
By the observation following Lemma \ref{tkmk}, 
the latter inequality implies that $t_k \ge m_k - n_k - C$, for some positive integer $C$ depending only on 
$\beta$. 
Consequently, we have
$$
\eps_{n_k + h} (x, \beta) = \eps_{ h} (x, \beta), \quad \hbox{for $h = 1, \ldots , \lfloor \delta n_k \rfloor$.} 
$$
This implies that the word $\eps_1(x, \beta) \ldots \eps_{n_k} (x, \beta)$ is the concatenation of prefixes 
of $\eps_1(x, \beta) \ldots \eps_{n_{k_0}} (x, \beta)$. 
By increasing $k_0$ if necessary, we can assume that the digits 
$\eps_1(x, \beta), \ldots , \eps_{\lfloor \delta n_{k_0} \rfloor } (x, \beta)$ are not all $0$. 
Let $L_k$ denote the length of the longest block composed solely of the digit $0$ in the word
$$
\eps_1 (x, \beta) \ldots \eps_{n_k} (x, \beta) \eps_1 (x, \beta) \ldots \eps_{n_k} (x, \beta).
$$
The sequence $(L_k)_{k \ge k_0}$ is bounded from above by $n_{k_0} + \lfloor \delta n_{k_0} \rfloor$. 
Consequently, the sequence $(m_k - t_k)_{k \ge 1}$ is also bounded from above. 
Set $g = \lceil \beta \rceil$ and consider the real number $\tx$ whose $g$-ary expansion 
is given by the infinite word $\eps(x, \beta)$. Since the sequence $(m_k - t_k)_{k \ge 1}$ is bounded, 
we have $\hr_\beta (x) = \hr_g (\tx)$, which is at most equal to $1$ by \eqref{rghat}. 
We have reached a contradiction. This proves the theorem.

\section{Proof of Theorem \ref{thdim}  \label{proofTh3}}

Let $r, \hr$ be real numbers such that $0 < r / (1 + r) < \hr \le 1$. 
Let $\eps$ in $(0, 1/3)$ be such that $0 < (r+ \eps) / (1 + r + \eps) < \hr - \eps$ and 
$$
\hr > 2 \eps > \eps, \quad r < 2 / \eps, \quad r > 5.
$$
We will show that the set of $x$ in $[0, 1)$ such that $r_\beta (x) \le r$ and $\hr_\beta (x) \le \hr$ has 
Hausdorff dimension zero. Let $x$ be a real number in this set. 
According to \cite[Lemma 3.3]{ZhWu20}, if $\eps$ is small enough, then we have
$$
n_{k+1} \le (1 - \eps) m_k, \quad (\hr - \eps) n_{k+1} \le m_k - n_k \le (r + \eps) n_k,
$$
for $k$ sufficiently large. Here and below, the sequences $(m_k)_{k \ge 1}, (n_k)_{k \ge 1}$, and 
$(t_k)_{k \ge 1}$ associated to $x$ are defined as in Section \ref{sec:3}. 

Assume first that $t_k \ge n_{k+1}$ for $k$ large enough (this is the case in particular 
when $(m_k - t_k)_{k \ge 1}$ is bounded). 
Then, we have
$$
\eps_{n_k+1} \ldots \eps_{n_k + \lfloor \eps n_k \rfloor} = \eps_{1} \ldots \eps_{ \lfloor \eps n_k \rfloor}, \quad 
\eps_{n_{k+1}+1} \ldots \eps_{n_{k+1} + \lfloor \eps n_{k+1} \rfloor} = \eps_{1} \ldots \eps_{ \lfloor \eps n_{k+1} \rfloor}.
$$
Assume that $n_{k+1} \le (1 + \eps (2r)^{-1} ) n_k$.  
Set $u = n_{k+1} - n_k$ and $h = \lfloor 2 r - 2 \rfloor$. Then, 
$$
\eps_1 \ldots \eps_u = \eps_{ju + 1} \ldots \eps_{(j+1)u}, \quad 1 \le j \le h-1.
$$
This implies that $r(x) \ge 2 r - 3 > r$, a contradiction. 
Consequently, there exists $C > 1$ such that $n_k \ge C^k$. 

Now, we use a covering argument to conclude that the Hausdorff dimension of the set 
of $x$ in $[0, 1)$ such that $r_\beta (x) \le r$ and $\hr_\beta (x) \le \hr$ is $0$. 
Assume that the prefix of length $n_1$ of $\eps(x, \beta)$ is given. Clearly, it yields the prefix of length $m_1$ 
and, since $n_{2} < m_1$, there are no more than
$m_1 - n_1$ choices for $n_2$. So, we get at most $(m_1 - n_1) \cdots (m_k - n_k)$ 
interval of length at most $\beta^{- m_{k+1}}$, thus at most $(2r)^k n_1 \ldots n_k$ 
intervals of length at most $\beta^{-C^k}$. 
Since $n_k \ge C^k$ for $k$ large enough, the series $\sum_{k \ge 1} (r n_k)^k \beta^{- s n_k}$ converges for every 
$s > 0$. The Cantelli lemma then implies that the Hausdorff dimension of our set is zero.

Let us explain now briefly why the same covering argument works 
as well in the case of an arbitrary sequence $(t_k)_{k \ge 1}$. 
The point is to show that $(n_k)_{k \ge 1}$ does not increase too slowly. 
If $t_k < n_{k+1}$ for many consecutive values of $k$ and 
$(n_k)_{k \ge 1}$ increases very slowly, then the infinite $\beta$-expansion of $1$ must be periodic 
and we further derive that the $\beta$-expansion of $x$ must be periodic, a contradiction. 
Consequently, $(n_k)_{k \ge 1}$ does not increase too slowly and the covering argument 
given above allows us to conclude in this case.

\section{Acknowledgements}
I started this work as I was enjoying the hospitality of NCTS Taiwan. 
I am very grateful to this Center for its support and the very good working conditions. 
I am also very thankful to the referees and the associate editor for their 
numerous suggestions and corrections.

\end{document}